\documentclass[12pt]{amsart}
\usepackage{amsmath}
\usepackage{amssymb}
\usepackage{amscd}
\usepackage{amsthm}

\numberwithin{equation}{section}
\theoremstyle{plain}
\newtheorem{theorem}[equation]{Theorem}

\newtheorem{proposition}[equation]{Proposition}
\newtheorem{lemma}[equation]{Lemma}

\theoremstyle{remark}

\theoremstyle{definition}
\newtheorem{definition}[equation]{Definition}

\newcommand{\ra}{\rightarrow}
\newcommand{\restr}{|}

\newcommand{\B}{{\mathcal B}}
\newcommand{\C}{\mathbb C}

\renewcommand{\H}{\mathbb H}

\newcommand{\N}{\mathbb N}

\newcommand{\R}{\mathbb R}

\renewcommand{\S}{{\mathbb S}}

\newcommand{\U}{{\mathcal U}}
\newcommand{\V}{{\mathcal V}}
\newcommand{\W}{{\mathcal W}}

\newcommand{\im}{\operatorname{Im}}
\newcommand{\ner}{\operatorname{Ner}}
\newcommand{\acts}{\curvearrowright}

\def\for{\quad \text{for}\quad}
\def\:{\colon}
\def\cangle{\tilde\angle}

\def\sub{\subseteq}
\newcommand{\al}{\alpha}
\def\de{\delta}
\def\De{\Delta}

\def\eps{\epsilon}

\def\la{\lambda}
\def\La{\Lambda}

\def\ra{\rightarrow}

\def\si{\sigma}

\def\geo{\partial_{\infty}}

\def\defeq{:=}

\def\Diam{\operatorname{diam}}
\def\diam{\operatorname{diam}}
\def\dist{\operatorname{dist}}
\def\Trip{\operatorname{Tri}}
\def\CAT{\operatorname{CAT}}
\def\WT{\operatorname{WT}}
\def\id{\operatorname{id}}
\def\Isom{\operatorname{Isom}}
\def\topdim{\operatorname{dim}_{\rm top}}
\renewcommand{\implies}{\Rightarrow}
\newcommand{\trip}{\operatorname{Tri}}

\newcommand{\ol}{\overline}
\newcommand{\no}{\noindent}
\newcommand{\qs}{quasi-sym\-metric}
\newcommand{\qm}{quasi-M\"obius}

\newcommand{\M}{M\"obius}

\begin{document}

\title{Rigidity for quasi-M\"obius  group actions}

\author{Mario Bonk}
\address{Department of Mathematics, University of Michigan,
Ann Arbor, MI, 48109-1109}
\email{mbonk@math.lsa.umich.edu}
\thanks{M.B.\ was supported by a Heisenberg fellowship
of the Deutsche Forschungsgemeinschaft.
B.K.\ was supported by NSF grant DMS-9972047}

\author{Bruce Kleiner}
\address{Department of Mathematics, University of Michigan,
Ann Arbor, MI, 48109-1109}
\email{bkleiner@math.lsa.umich.edu}

\begin{abstract}
If a group acts  by uniformly \qm\ homeomorphisms on a compact
 Ahlfors $n$-regular space of topological dimension $n$ such that
the induced action on the space of distinct triples is cocompact, then the 
action is quasi-symmetrically conjugate to an action
on the standard $n$-sphere by M\"obius transformations.   
\end{abstract}

\date{December 14, 2001}
\maketitle


\section{Introduction}

It has been known since the time of Poincar\'e that the limit set 
of a subgroup of ${\rm PSL}(2,\C)$ obtained by a small 
deformation of a discrete cocompact subgroup of
${\rm PSL}(2, \R)\sub{\rm PSL}(2,\C)$
will be a nowhere differentiable curve unless it is round. 
Much later   R.\ Bowen \cite{bowen} made this more precise 
by proving 
that such a limit curve is either  a round circle or has Hausdorff dimension
strictly greater than 1.
The group 
 ${\rm PSL}(2,\C)$ is  isomorphic to the group of orientation preserving 
isometries of hyperbolic $3$-space. Therefore, it is a 
natural question whether similar results hold for subgroups 
of the isometry group $\Isom(\H^{n+1})$ of hyperbolic 
$(n+1)$-space when   $n\ge 2$, or, what is the same, for groups 
of  M\"obius transformations acting on the standard 
$n$-sphere $\S^n$.  
Rigidity results in this vein
 were obtained by
Sullivan \cite[p.\ 69]{sul1} and 
Yue \cite[Theorem 1.5]{yue}. 

In the present paper we generalize these results further by considering
uniformly \qm\  group actions on compact metric spaces $Z$ that 
induce cocompact actions  
on the space $\trip(Z)$ of distinct triples of $Z$.  
The following theorem is our  main result.

\begin{theorem}
\label{mainthm}
Let $n\in\N$, and let $Z$ be a compact, 
Ahlfors $n$-regular metric space
of topological dimension $n$.  Suppose  $G\acts Z$  
 is a uniformly \qm\  action of a group $G$ 
on $Z$, where the induced action  $G\acts \Trip(Z)$
is  cocompact. 
Then   $G\acts Z$
is quasi-symmetrically conjugate to an  action
of $G$ on the standard sphere $\S^n$ by \M\ transformations.
\end{theorem}

The terminology will be explained in the body of the paper. 
Note that part of the conclusion is that
 $Z$ is  homeomorphic to $\S^n$.

When $G$ is a hyperbolic group, the boundary
$\geo G$ carries a metric $d$  unique
up to quasi-symmetry, with respect to which the 
canonical action $G\acts\geo G$ is 
uniformly \qm.
In this case the induced action on $\Trip(\geo G)$ 
is discrete and cocompact, so Theorem \ref{mainthm}
may be applied if $(\geo G,d)$ is \qs\ to an
Ahlfors $n$-regular space whose topological 
dimension  is equal to $n$.
Note that $(\geo G,d)$ will always be Ahlfors $Q$-regular for some 
$Q>0$, but in general $Q$ will exceed the topological dimension 
of $\geo G$.

In order to state our next result, 
we recall (see the discussion in Section \ref{hypgp}) that if $X$ is a
$\CAT(-1)$-space,
then any point $p\in X$ determines a canonical   metric on
$\geo X$, and any two such metrics are bi-Lipschitz
equivalent by the identity map.   In particular,
we may speak  of the Hausdorff dimension of any subset of  $\geo X$, 
since this number is independent of the choice of the 
canonical metric. 
We then have the following corollary of Theorem~\ref{mainthm}  which 
generalizes  a result by Bourdon
\cite[0.3 Th\'eor\`eme ($\H^n$ case)]{Bou2}.

\begin{theorem}
\label{cat-1theorem}
Suppose $n\in \N$, $n\ge 2$. 
Let $G\acts X$ be a properly discontinuous, quasi-convex cocompact,
and isometric
action on a $\CAT(-1)$-space $X$.  If the Hausdorff dimension
and topological dimension of the limit set $\Lambda(G) \sub
 \geo X$ are both equal to 
$n$,  then $X$ contains a convex, $G$-invariant subset  $Y$
isometric to $\H^{n+1}$ on which $G$ acts 
 cocompactly.
\end{theorem}

The terminology and the notation will be explained in 
Section~\ref{hypgp}. 
Note that  the ineffective
kernel $N$ of the induced action $G\acts Y$ is finite,
and $G/N$ is isomorphic to a uniform lattice in $\Isom(\H^{n+1})$.

In contrast to Theorem \ref{mainthm} where the case $n=1$ is allowed, we 
assume $n>1$ in the previous theorem, in order to be able to apply
Bourdon's result. It is an interesting question whether the statement is 
also true for $n=1$. 
See  the discussion in
Section~\ref{hypgp}.

The proof of  Theorem \ref{mainthm} can be outlined as follows. First, 
we use the dimension assumption to get a Lipschitz map 
$f\:Z\ra\S^n$ such 
that the image   of $f$  has positive Lebesgue measure.
According to a result by 
David and Semmes 
one can rescale $f$ and extract
a limit mapping $\phi\:X\ra\R^n$ defined on a weak tangent space of 
$Z$ which has
bounded multiplicity, i.e.\ point
inverses $\phi^{-1}(y)$ have uniformly
bounded cardinality.   We then show that $\phi$
is locally bi-Lipschitz somewhere, and
as a consequence some weak tangent of $Z$
is bi-Lipschitz to $\R^n$. The  assumptions on the group action 
can then be used to prove that 
$Z$ is \qs\ to $\S^n$. Once this is established, the theorem  follows 
from a result by Tukia.

Our method of  proving  Theorem~\ref{mainthm} can also be applied in
other contexts. In \cite[Question 5]{heinsem}
Heinonen and Semmes  ask whether
every linearly locally contractible  Ahlfors $n$-regular 
metric $n$-sphere $Z$  that is quasi-symmetrically three point homogeneous 
is quasi-symmetrically equivalent to the standard $n$-sphere $\S^n$.
One can show that  the answer to this question is 
positive,
 if we make the stronger assumption that  $Z$ is 
 three point homogeneous by uniform \qm\ homeomorphisms. 
(see the discussion in  Section~\ref{proofmainthm}).

\medskip \noindent 
{\bf Acknowledgement.}
A previous version of this paper was  based on some 
rather deep results on the  uniform rectifiability 
of metric spaces satisfying some topological nondegeneracy assumptions. 
The statements   we needed are  implicitly contained in the works of David 
and Semmes, but not stated  explicitly.   
The approach taken in this version uses a much more elementary result by
David and Semmes. 
The  authors are indebted
to Stephen Semmes for conversations about these issues and 
thank him especially for directing their attention to the results in Chapter 12 
of \cite{davsem}.

\medskip \noindent
{\bf Notation.}
The following  notation will be used throughout the paper. 
Let $Z$ be a metric space. The metric on $Z$ will be denoted 
by $d_Z$, and the open and the closed ball of radius $r>0$
centered at $a\in Z$ by $B_Z(a,r)$ and $\bar B_Z(a,r)$,
respectively. We will drop the subscript $Z$ if the space $Z$ is understood.
If $A\sub Z$ and $d=d_Z$, then $d\restr_A$ is the restriction of the metric
$d$ to $A$. 
 We use 
$\diam(A)$ for the diameter, $\bar A$ for the closure, and 
$\#A$ for the cardinality of a set $A\sub Z$.  
If $z\in Z$ and  $A,B\sub Z$, then  
$\dist(z,A)$ and   $\dist(A,B)$ are  the  distances of $z$ and $A$ and
of  $A$ and $B$, respectively.  
If $A\sub Z$ and $r>0$, then we  let  $N_r(A):=\{z\in Z: \dist(z, A)<r\}$. 
The
{\em Hausdorff distance} of two sets $A,B\sub Z$ is defined by
$$ \dist_H(A,B):= \max\big\{ \sup_{a\in A} \dist(a, B),\
\sup_{b\in A} \dist(b, A)\big\}. $$

Suppose $X$ and $Y$ are metric spaces. If $f\:X\to Y$ is a map, then 
we let $\im(f):=\{f(x): x\in X\}$. If $A\sub X$, then $f\restr_A$ denotes
the restriction of the map $f$ to $A$.
If $g\:X\to Y$ is another map, we
 let
$$ \dist(f,g):=\sup_{x\in X} \dist(f(x), g(x)).$$
The identity map on a set  $X$ will be denoted by $\id_X$.

\section{Quasi-M\"obius maps and group  actions}
\label{qm}

Let $(Z,d)$ be a metric space.
   The {\em cross-ratio}
of a four-tuple of  distinct points $(z_1,z_2,z_3,z_4)$ in $Z$
is the quantity
$$[z_1,z_2,z_3,z_4]:= \frac{d(z_1,z_3)d(z_2,z_4)}{d(z_1,z_4)d(z_2,z_3)}.$$

Suppose $X$ and $Y$ are metric spaces.
Suppose  $\eta\:[0,\infty)\ra[0,\infty)$
is a homeomorphism, and let   $f\:X\ra Y$ 
be an   injective map.
The map $f$ is an {\em $\eta$-quasi-\M\ map} if for
every four-tuple $(x_1,x_2,x_3,x_4)$ of distinct
points in $X$, we have
$$[f(x_1),f(x_2),f(x_3),f(x_4)]\leq \eta([x_1,x_2,x_3,x_4]).$$
Note that by exchanging the roles of $x_1$ and $x_2$, one
gets the lower bound
$$\eta([x_1,x_2,x_3,x_4]^{-1})^{-1}\leq [f(x_1),f(x_2),f(x_3),f(x_4)].$$
Hence $f$ is a homeomorphism onto its image $f(X)$,
and the inverse map $f^{-1}\:f(X)\ra X$ is also \qm.

The map $f$ is  {\em $\eta$-\qs\ } if
$$\frac{d_Y(f(x_1),f(x_2))}{d_Y(f(x_1),f(x_3))}
\leq \eta\left (\frac{d_X(x_1,x_2)}{d_X(x_1,x_3)}\right )$$
for every triple $(x_1,x_2,x_3)$ of distinct points in $X$.

Finally, $f$ is called {\em bi-Lipschitz} if there exists a constant $L\ge 1$
(the bi-Lipschitz constant of $f$) such that 
$$ (1/L) d_X(x_1,x_2) \le d_Y(f(x_1), f(x_2))\le L d_X(x_1,x_2),$$
whenever $x_1,x_2\in X$.

We mention some basic properties of these maps. 

(1) The post-composition of an $\eta_1$-\qm\ map by  an
$\eta_2$-\qm\ map is an $\eta_2\circ\eta_1$-\qm\ map.
Similar statements are true for \qs\ maps and bi-Lipschitz maps. 

(2) A bi-Lipschitz map is \qs\ and \qm. A \qs\ map is \qm. 
 A \qm\ map defined on a bounded space is \qs.

(3) Let $X$ and $Y$ be  compact metric spaces, and suppose
  $f_k\:X\ra Y$ is an   $\eta$-\qm\ 
map for $k\in \N$. Then we have that 

 (a) the sequence $(f_k)$ subconverges uniformly to an $\eta$-\qm\
map,  or

(b) there is a point $x_0\in X$ 
so that the sequence  $(f_k\restr_{X\setminus\{x_0\}})$
subconverges uniformly on 
compact subsets of $X\setminus\{x_0\}$  to a constant map.

The  alternative (b) can be excluded by a normalization 
condition; namely, that each map $f_k$ maps a uniformly
separated triple of points in $X$ to a uniformly
separated triple in $Y$.

\medskip
We will need the following  extension of property (3).

%
%
%

\begin{lemma}
\label{TBA}
Suppose $(X, d_X)$ and $(Y,d_Y)$ are  compact
 metric spaces, and let
$f_k\:D_k\ra Y$ for $k\in \N$ be an $\eta$-\qm\ map defined on
a subset $D_k$ of $X$. Suppose
$$ \lim_{k\to \infty}\,  \dist_H(D_k, X) = 0,$$
and that  for $k\in \N$  
there exist  triples $(x^1_k, x^2_k,x^3_k)$
and $(y^1_k, y^2_k,y^3_k)$ of points in $D_k$
and $Y$, respectively,
such that
$$f_k(x^i_k)=y^i_k \for  k\in \N, \ i\in \{1,2,3\},$$ 
$$ d_X(x^i_k,x_k^j)\ge \de \ \text{and}\
d_Y(y^i_k,y_k^j)\ge \delta 
\text{ for } k\in \N, \, i,j\in \{1,2,3\},\, i\ne j, $$
where $\delta>0$ is independent of $k$. 

Then the sequence $(f_k)$ subconverges uniformly
 to a \qm\ map $f\:X\ra Y$, i.e.\ there exists a
monotonic sequence $(k_\nu)$ in $\N$  such that
$$ \lim_{\nu\to \infty} \dist(f_{k_\nu}, f\restr_{D_{k_\nu}})=0. 
$$
Suppose in addition that 
$$ \lim_{k\to \infty}\,  \dist_H(f_k(D_k), Y) = 0.$$
Then the sequence  $(f_k)$ subconverges uniformly
 to  \qm\ homeomorphism $f\:X\ra Y$.
\end{lemma}

The lemma says that a sequence  $(f_k)$ 
of uniformly \qm\ maps defined   on denser and denser  
 subsets of a space $X$  and mapping into the 
same space $Y$ subconverges to a \qm\ map defined on the whole space
$X$, if each map $f_k$ maps a uniformly separated triple in $X$ to 
a uniformly separated triple in $Y$.
Moreover,  a surjective  limiting map can be obtained if 
the images of the maps $f_k$ Hausdorff converge to  the space  $Y$.

\proof 
The assumptions imply that the functions  $f_k$ are equicontinuous (cf.\
\cite[Thm.\ 2.1]{Vai}).  The proof of the first part of the  lemma
 then follows from standard arguments based on the Arzel\`a-Ascoli 
theorem, and we leave the details  to the reader.

To prove the second part, note that according 
to the first part, by passing to a subsequence if necessary,  
we may  assume that 
$$ \dist(f_k, f\restr_{D_k})\to 0 \for k\to \infty.$$ 
Let $D_k':=f_k(D_k)$ and $g_k:=f_k^{-1}\:D_k'\to X$. 
The maps $g_k$ are uniformly \qm. Hence,  by our  additional 
assumption we can apply the first part of the lemma to the sequence
$(g_k)$.
Again by selecting a  subsequence  of $(g_k)$  if necessary,
we may assume  that 
$$ \dist(g_k, g\restr_{D'_k})\to 0 \for k\to \infty,$$ 
where $g\: Y\to X$ is a \qm\ map. 
Since $g_k\circ f_k =\id_{D_k}$ and $f_k\circ g_k =\id_{D'_k}$, we obtain 
from the uniform convergence of the sequences $(f_k)$ and $(g_k)$ 
that $g\circ f =\id_{X}$ and $g\circ f=\id_{Y}$. Hence $f$ is a bijection 
and therefore a \qm\ homeomorphism. \qed
\medskip

Let $Z$ be an unbounded locally compact metric space with metric
$d=d_Z$, let  $p\in Z$
be a base point, and let
 $\hat Z=Z \cup\{\infty\}$ be  the one-point compactification
of $Z$. In order to define a metric on $\hat Z$ associated with the 
pointed space $(Z,p)$ let $h_p\: \hat Z \to [0,\infty)$ be given by 
\begin{equation*}
h_p(z): = \left\{ 
\begin{array}{ll} \displaystyle\frac 1{1+d(z,p)}& \for z\in Z,  \\
\quad\quad0& \for z=\infty.  
\end{array} \right. 
\end{equation*}

Moreover, let 
$$ \rho_p(x,y) = h_p(x)h_p(y) d(x,y)\for x,y\in Z, $$
$\rho_p(x,\infty)=\rho_p(\infty, x)= h_p(x)$
for $x\in Z$, $\rho_p(\infty,\infty)=0$.
Note that if an argument of the functions  $h_p$ and $\rho_p$ is the point
at infinity, the corresponding value  can 
be obtained  
as a limiting case of
values at arguments in $Z$.
Essentially, 
the function $\rho_p$ is the metric on $\hat Z$ that we are
looking for. This distance function is an analog 
of the chordal metric on the Riemann sphere. 
 Unfortunately, $\rho_p$ will not satisfy the triangle inequality 
in general. We remedy this problem by a standard procedure. 

If $x,y\in \hat Z$ we define 
$$ \hat d_p(x,y):= \inf \sum_{i=0}^{k-1} \rho_p(x_i, x_{i+1}), $$
where the infimum is taken over all finite sequence of points $x_0, \dots,
x_k\in \hat Z$ with $x_0=x$ and $x_k=y$. 

\begin{lemma} \label{dhat}
 The function $\hat d_p$ is a metric on   $\hat Z$ whose 
induced topology agrees with the topology of $\hat Z$. 
The identity map $\id_Z\: (Z,d)\to (Z, \hat d_p\restr_Z) $ is an $\eta$-\qm\ 
homeomorphism where $\eta(t)=16 t$.
\end{lemma}

\proof The first part of the lemma  immediately  follows if we can show that
\begin{equation} \label{metricAC}
 \frac  14  \rho_p(x,y) \le \hat d_p(x,y) \le \rho_p(x,y) \for
x,y\in \hat Z.
\end{equation} 
The second part also follows from this inequality  by observing that 
if $(z_1, z_2, z_3, z_4)$ is a four-tuple of distinct points 
in $Z$, then 
$$ 
 \frac{\hat d_p(z_1,z_3)\hat d_p(z_2,z_4)}{\hat d_p(z_1,z_4)\hat d_p(z_2,z_3)}
\le 16 \frac{\rho_p(z_1,z_3)\rho_p(z_2,z_4)}{\rho_p(z_1,z_4)\rho_p(z_2,z_3)}
= 16 \frac{d(z_1,z_3)d(z_2,z_4)}{d(z_1,z_4)d(z_2,z_3)}.
$$

The right  hand inequality in \eqref{metricAC}
 follows from the 
definition of $\hat d_p$. In order to prove the left hand inequality, 
we may assume $h_p(x)\ge h_p(y)$ without loss of generality. 
Moreover, we may assume $x\in Z$ and so  $h_p(x)>0$, because otherwise
$x=y=\infty$ and the inequality is true.

If $x_0, \dots, x_k$ is an arbitrary sequence with $x_0=x$ and $x_k=y$, 
we consider two cases:

If $h_p(x_i)\ge \frac 12 h_p(x)>0$ for all $i\in \{0, \dots, k\}$,
then $x_i\in Z$,
and the triangle inequality applied to $d$ gives

\begin{eqnarray}\label {lbdsum1}
\sum_{i=0}^{k-1} \rho_p(x_i, x_{i+1})& \ge&
\frac 14 h_p(x)^2\sum_{i=0}^{k-1} d(x_i, x_{i+1}) \\ 
&\ge&  \frac 14 d(x,y) h_p(x) h_p(y)
=\frac 14 \rho_p(x,y). \nonumber
\end{eqnarray}

Suppose there exists $j\in \{0, \dots, k\}$ such that 
$h_p(x_j) <\frac 12 h_p(x)$. Note that that it follows from the definitions
that   $|h_p(u)-h_p(v)| \le \rho_p(u,v)$ for $u,v\in \hat Z$. 
Moreover, since $h_p(y) \le h_p(x) $ we have 
$d(x,p)\le d(y,p)$ in case $y\in  Z$. This implies 
$$ \frac {d(x,y)} {1+d(y,p)} \le 2 \frac {d(y,p)} {1+d(y,p)} \le 2,$$
which leads to  
$\rho_p(x,y) \le 2h_p(x)$. This is also true if $y=\infty$. 
We arrive at 
\begin{equation}\label {lbdsum2}
\sum_{i=0}^{k-1} \rho_p(x_i, x_{i+1}) \ge \sum_{i=0}^{k-1}
|h_p(x_i)-h_p(x_{i+1})|\ge \frac 12 h_p(x) \ge  \frac 14 \rho_p(x,y).
\end{equation} 
The desired inequality follows from \eqref{lbdsum1} and \eqref{lbdsum2}. 
\qed 

\medskip

Let $(Z,d)$ be a metric space. We write $G\acts Z$, if 
$G$ is a group that acts on $Z$ by homeomorphisms.
The image of a point $z\in Z$ under the  group element $g$ is denoted
by $g(z)$. 
The action $G\acts Z$ is called  {\em faithful}
if the only element in $G$ that  acts as the identity on $Z$ is the unit
element.

If $\eta\:[0,\infty)\ra[0,\infty)$
is a homeomorphism, then an action $G\acts Z$
is an {\em $\eta$-\qm\ action} if each $g\in G$
induces an $\eta$-\qm\ homeomorphism of $Z$.
An action $G\acts Z$ is {\em uniformly \qm\ }
if it is $\eta$-\qm\ for some homeomorphism
$\eta\:[0,\infty)\ra[0,\infty)$.
If $Z$ is locally compact, then the action $G\acts Z$ is called 
{\em cocompact} if there exists a compact set $K\sub Z$ such that 
$$Z=\bigcup_{g\in G} g(K).$$ 

We denote  by 
$$\trip(Z):= \{ (z_1,z_2,z_3)\in Z^3: z_1\ne z_2\ne z_3\ne z_1\}$$
the space of distinct triples in $Z$. If $G\acts Z$ is a group action,
then there is a natural induced action  $ G\acts \trip(Z)$
defined by 
$$ g(z_1,z_2,z_3) := (g(z_1), g(z_2), g(z_3)) $$ 
for $g\in G$ and $(z_1,z_2,z_3)\in \trip(Z)$. 

Suppose $G \acts Z$ is an action on a compact space $Z$. Then
the induced  action $ G\acts \trip(Z)$ is cocompact if and only if there 
exists $\delta>0$ such that 
for every triple $ (z_1,z_2,z_3)\in \trip(Z)$ there exists a group element  
 $g\in G$ such that
$$ d(g(z_i), g(z_j)) \ge 
 \delta \for i,j\in \{1,2,3\},\ i\ne j.$$ 
This condition means that every triple in $\trip(Z)$ can be mapped 
to a uniformly separated triple by  some map $g\in G$.

%

\section{Maps of bounded multiplicity}

The goal of this section is to study
continuous maps of bounded multiplicity  between 
a space of topological dimension $n$  and $\R^n$.
The main result is
Theorem~\ref{covertheorem} which may be of independent interest.

\begin{definition}
If $f\:X\ra Y$ is a continuous  map between  metric
spaces  $X$ and $Y$, then $y\in Y$ is a {\em stable value}
of $f$ if there is  $\eps>0$ such that
$y\in \im(g)$ for every continuous   map $g\:X\ra Y$ with $\dist(f,g)<\eps$.
\end{definition}
Note that  
 the set of stable values of a map $f\:X\ra\R^n$ is
an open subset of $\R^n$.

Recall that a map is {\em light} if all point
inverses are totally disconnected.
We will prove the following proposition. 
\begin{proposition}
\label{stablevalue}
Let $X$ be a   compact
metric space of topological dimension at least  $n$,
and let $f\:X\ra \R^n$ be a light continuous   map.
Then $f$ has stable values.
\end{proposition}

As we will see,  the proof is a slight amplification of the 
well-known argument that such a map $f$ cannot decrease
topological dimension.

\begin{definition}
A map $f\:X\ra Y$ between two spaces has {\em bounded
multiplicity} if there is a constant
$N\in\N$ such that $\#f^{-1}(y)\leq N$ for all $y\in Y$.
\end{definition}

Using Proposition \ref{stablevalue}
we will prove:

\begin{theorem}
\label{covertheorem}
Suppose $X$ is a compact metric space, 
every nonempty  open subset of $X$ has  topological dimension at least
$n$,  and $f\:X\ra\R^n$ is a
continuous map of  bounded multiplicity. 
Then there is an open subset $V\sub \im(f)$
with $\bar V=\im(f)$,
such that $U:=f^{-1}(V)$
is dense in $X$ and $f\restr_U\:U\ra V$
is a covering map. 
\end{theorem}

In particular, there exist nonempty open sets $U_1\sub X$ and $V_1\sub \R^n$
such that $f\restr_{U_1}$ is a homeomorphism of $U_1$ onto $V_1$. 
It is in this form that we will use Theorem~\ref{covertheorem} in the proof 
of Theorem~\ref{mainthm}.

\medskip

Let $X$ be a topological space, and let 
${\mathcal U}=\{U_i: i\in I\}$ be  a  cover of $X$ by open
 subsets
$U_i$  indexed by  some set $I$. 
The {\em nerve} of ${\mathcal U}$,  denoted by $\ner({\mathcal U})$,
is a simplicial complex whose simplices corresponds
to the subsets $I'\sub I$ for which
$$ U_{I'}:= \bigcap_{i\in I'} U_i \ne \emptyset. $$
The {\em order} of ${\mathcal U}$ is the supremum of all  numbers
$\#I'$ such that  $U_{I'}\ne \emptyset$.
 We denote the topological dimension 
of $X$ by $\topdim(X)$ (cf.\ \cite[Def.~I.4]{nag}).   
A compact metric space  $X$ has  topological dimension at most $n$, 
if and only if 
 open covers of order at most $n+1$ are cofinal in the family
of all open covers of $X$, i.e.,  every
open cover has an  open refinement which has order at most $n+1$. 
The order of an open cover ${\mathcal U}$
   is equal to $\topdim(\ner({\mathcal U}))+1$.

In order to prove Proposition~\ref{stablevalue} we discuss a general 
construction that associates  a fine cover with a  light continuous 
 map
$f\:X\ra Y$  from a  compact metric space
$X$ to a separable metric space $Y$.  Pick $\eps>0$. 

If $y\in Y$, then $f^{-1}(y)$ is compact and totally disconnected,
so the diameter of connected components of $N_\de(f^{-1}(y))$
tends to zero as $\de\ra 0$.  Hence there is a number $r_y>0$
such that $N_{r_y}(f^{-1}(y))$ can be decomposed as a finite
disjoint union of open sets with diameter  less than $\eps$;
moreover, there is a number $s_y>0$ such that $f^{-1}(B(y,s_y))\sub
N_{r_y}(f^{-1}(y))$.  Let $\B$ be a finite 
cover of $\im(f)$ by  balls of the form $B(y,s_y)$.

Suppose $\U=\{U_i: i\in I\}$
is a cover of $\im(f)$ by open subsets of $Y$. Let  $I':=\{ i\in I: U_i\cap 
\im(f)\ne \emptyset\}$, and assume that  $\U':=\{ U_i: i\in I'\}$ 
 refines  $\B$.
Then $f^{-1}(\U'):=\{f^{-1}(U_i): i\in I'\}$ is an
open cover of $X$ such that for all $i\in I'$,
we have $f^{-1}(U_i)\sub N_{r_y}(f^{-1}(y))$ for 
some $y\in Y$, which implies that $f^{-1}(U_i)$
may be written as a finite disjoint union of open
subsets with diameter less than  $\eps$.  Choosing such a
decomposition of $f^{-1}(U_i)$ for each $i\in I'$ yields
a  collection of open sets $\V=\{V_j: j\in J\}$
which covers $X$, and a map $\al\:J\ra I'\sub  I$  such that 
$V_j$ is an open set appearing in the decomposition of
$f^{-1}(U_{\al(j)})$.
 Note that $\al$ induces a simplicial
map $\phi:\ner(\V)\ra \ner(\U)$ since

$$V_{j_1}\cap\ldots\cap V_{j_k}\neq\emptyset\implies
f^{-1}(U_{\al(j_1)})\cap\ldots\cap f^{-1}(U_{\al(j_k)})\neq\emptyset$$
$$\implies U_{\al(j_1)}\cap\ldots\cap U_{\al(j_k)}\neq\emptyset.$$
In fact, $\phi$ is injective on simplices, since if
$j,\,j'\in J$ are distinct and
$\al(j)=\al(j')$,
then $V_j$ and $V_{j'}$ are disjoint fragments of the
same open set $f^{-1}(U_{\al(j)})=f^{-1}(U_{\al(j')})$,
and so $V_{j}\cap V_{j'}=\emptyset$.  In particular,
we have $\topdim(\ner(\V))\leq \topdim(\ner(\U))$.

Suppose $\{\rho_i: i\in I\}$ is a partition of unity in $Y$
 subordinate to 
$\U$. Here and in the following we interpret subordination 
in the sense that  $\{\rho_i \ne 0\}\sub U_i$ for all $i\in I$. 
We can produce a partition of unity $\{\nu_j: j\in J\}$ in $X$ 
subordinate to $\V$ as follows: let 
$\nu_j\defeq \chi_{{}_{\scriptstyle V_j}}
\cdot\left(\rho_{\al(j)}\circ f\right)$,
where $\chi_{{}_{\scriptstyle V_j}}$ is the characteristic function of $V_j$.
Using the functions  $\{\rho_i: i\in I\}$ 
as barycentric coordinates in $\ner(\U)$, and
the functions $\{\nu_j: j\in J\}$ as barycentric coordinates in 
$\ner(\V)$, we  obtain induced continuous  maps $\rho\: Y\to \ner(\U)$ and 
$\nu\: X \to \ner(\V)$ 
such that 
$\phi\circ\nu=\rho\circ f$.  

We note that since $\eps>0$
was chosen arbitrarily, if we have a cofinal
family  of covers $\U'$ of $\im(f)$ of order
at most $N$, then the corresponding family of covers 
$\V$ of $X$ will be cofinal and
its members will  have order at most  $N$; this implies
that $\topdim(Y) \ge \topdim(\im(f)) \ge \topdim(X)$.

\medskip
\noindent{\em Proof of Proposition \ref{stablevalue}.}
For $\eps>0$ we now apply the construction above
in the special case that $Y=\R^n$, $\topdim(X)\ge n$,
and the open cover 
$\U=\{U_i:i\in I\}$ of $\R^n$ is the open star cover associated with 
a  triangulation of $\R^n$. Since  $f(X)$ is compact, the associated cover 
$\U'$ will refine a given cover of $f(X)$ if the triangulation
of $\R^n$ is chosen fine enough.
We have a homeomorphism
$\rho\:\R^n\ra\ner(\U)$ (we conflate simplicial complexes with their
geometric realizations), and an induced partition of unity
$\{\rho_i: i\in I\}$ coming from the barycentric coordinate
functions of the map $\rho$.

Since the family of  open  covers of $X$ induced by our construction 
 is cofinal in the 
family of all open covers   of $X$, we can choose 
$\eps>0$  
small enough so that the induced cover  $\V$ of $X$ 
 does not admit an   open  refinement
$\W$ of order at most $n$. 

\begin{lemma}
Some $n$-simplex $\si$ of $\ner(\V)$ has an interior
point $\xi$ which is a stable value of $\nu\:X\ra\ner(\V)$.
\end{lemma}
\proof
Suppose not.  Then we may form  a set $S$ by choosing
one interior point from each $n$-simplex of $\ner(\V)$,
 and perturb $\nu$
slightly  on a small neighborhood of $\nu^{-1}(S)$
to get a map $\nu'\:X\ra\ner(\V)$  such that its barycentric
coordinate functions are subordinate to $\V$, and $\im(\nu')\cap
S=\emptyset$. (See the first part of the proof of
Lemma~\ref{existsgreen} for the idea of how to construct this perturbation.) 
 Then we may compose $\nu'$ with 
the ``radial projection" in each $n$-simplex to get a 
map $\nu''$ that maps $\ner(\V)\setminus S$ to
the $(n-1)$-skeleton $[\ner(\V)]_{n-1}$ of $\ner(\V)$
and  whose
barycentric coordinates are subordinate to $\V$;
pulling back the open star cover of $\ner(\V)$
by $\nu''$, we get a refinement of $\V$ of order
at most $n$, which is a contradiction.
\qed

\medskip
If $\xi$ is as in the lemma, then $\phi(\xi)\in \ner(\U)$
is clearly a stable value of $\phi\circ\nu\:X\ra\ner(\U)$;
but $f=\rho^{-1}\circ\phi\circ\nu$ where $\rho^{-1}$
is a homeomorphism, so $\rho^{-1}(\phi(\xi))$ is a  stable
value of $f$.  This completes  the proof of Proposition
\ref{stablevalue}. \qed

\begin{definition}
Let $X$ be a topological  space, and $f\:X\ra\R^n$ be a map.
Then $x\in X$ is a {\em stable point}  of $f$ if $f(x)$
is a stable value of $f\restr_U$ for every neighborhood
$U$ of $x$.  
\end{definition}

\begin{lemma}
\label{existsgreen}
Suppose $X$ is metric space, and $f\:X\ra\R^n$
is a continuous  map.  Then $y\in \R^n$ is a stable value of 
$f$ if and only if $y$ is a stable value of $f\restr_{f^{-1}(W)}$
for every  neighborhood $W$ of $y$.
When $X$ is a compact metric space and $f^{-1}(y)$
is totally disconnected, then $y$ is a
stable value of $f$ if and only if the fiber $f^{-1}(y)$
contains a stable point.
\end{lemma}

\proof
We will only prove the ``only if'' implications;
the other implications are immediate.

Suppose $W\sub\R^n$ is an open neighborhood
of $y$, and  $y$ is  an unstable value
of $f\restr_{U}$, where  $U:= f^{-1}(W)$.
Choose $\de>0$ such that $\bar B(y,\de)\sub W$,
and let $V\defeq f^{-1}(\R^n\setminus \bar B(y,\de))$.
Pick $\eps>0$.  As $y$ is  an unstable value of $f\restr_U$,
we can find a map $g_U\:U\ra\R^n$ such that
$\dist (g_U,f\restr_U)<\min(\eps,\de)$ and $y\not\in\im(g_U)$.
Define $g_V\:V\ra\R^n$ to be the restriction of 
$f$ to $V$.  Combining $g_U$ and $g_V$ using a partition of 
unity subordinate to the cover $\{U,V\}$, 
we get a continuous  map $g\:X\ra\R^n$ such that 
$\dist(g,f)<\eps$ and $g^{-1}(y)=\emptyset$.
Since $\eps>0$ was arbitrary, we have shown
that $y$ is  not a stable value of $f$.

Now suppose $X$ is compact, $f^{-1}(y)$ is totally
disconnected, and every point $x\in f^{-1}(y)$
is unstable.  By the compactness of 
$f^{-1}(y)$ we can find a finite cover 
$\B=\{B(x_1,r_1),\ldots,B(x_k,r_k)\}$ of $f^{-1}(y)$ by balls
 where $x_i\in f^{-1}(y)$
and $y$ is an unstable value of $f\restr_{B(x_i,r_i)}$
for each $1\leq i\leq k$. When $\de>0$ is
sufficiently small, then $f^{-1}(B(y,\de))$ can be
decomposed into a disjoint union of open sets
$U_1,\ldots,U_j$ so that the cover $\{U_i\}$
of $f^{-1}(y)$ refines $\B$.  This means that $y$ is 
an unstable value of $f\restr_{U_i}$ for 
each $i$, which implies that 
$y$ is an unstable value of $f\restr_{f^{-1}(B(y,\de))}$.
This is a  contradiction to what we proved in the first part  of the proof.
\qed

\medskip

Now let  $X$  be a compact metric
space such that $\topdim(U) \ge n$ for all nonempty open subsets
$U\sub X$, and $f\:X\ra \R^n$ be a  continuous map 
of bounded multiplicity.

\begin{lemma}
\label{semicontinuity}
For all $y\in\R^n$ and all $\eps>0$, there
is  $\de>0$ such that for all $y'\in B(y,\de)$
and all stable points $x\in f^{-1}(y)$, there
is a stable point in $f^{-1}(y')\cap B(x,\eps)$. 
\end{lemma}
\proof Let $ \{x_1, \dots, x_k\}$ be the stable points in $f^{-1}(y)$ and 
pick $i\in\{1,\ldots,k\}$.  Since $x_i$ is stable point,
$y$ is a stable value of $f\restr_{B(x_i,\eps)}$.
So any $y'$ sufficiently close to $y$ is also
a stable value of $f\restr_{B(x_i,\eps)}$ and by
Lemma \ref{existsgreen} for such $y'$ we will have a stable point in 
$f^{-1}(y')\cap B(x_i,\eps)$.  This holds for all $i$, so
the lemma follows.
\qed

\medskip

We define the {\em stable multiplicity function}
$\mu\:\R^n\ra\N$ by letting
$\mu(y)$ be the number of stable points in 
$f^{-1}(y)$.

\begin{lemma}
\label{allgreen}
If $\mu$ is locally maximal at $y\in\R^n$,
then every $x\in f^{-1}(y)$ is stable.
\end{lemma}
\proof
Let $U\sub \R^n$ be a neighborhood of $y$ such 
that $\mu(y')\leq \mu(y)$ for all $y'\in U$.
Let $x_1,\ldots,x_k$ be the stable points in 
$f^{-1}(y)$, and suppose 
$x\in f^{-1}(y)\setminus\{x_1,\ldots,x_k\}$.
Pick $\eps>0$ such that the balls
$B(x,\eps),\,B(x_1,\eps),\ldots,B(x_k,\eps)$
are disjoint. 

Choose $\de>0$ as in the previous lemma. 
Let 
$y'$ be a stable value of $f\restr_{B(x,\eps)}$
lying in $U\cap B(y,\de)$; such a $y'$ exists since
by Proposition \ref{stablevalue} stable values of $f\restr_{B(x,\eps)}$
are dense in $\im(f\restr_{B(x,\eps)})$.
Then $f^{-1}(y')$ has a stable point in each of the
balls $B(x,\eps),\,B(x_1,\eps),\ldots,B(x_k,\eps)$,
so $\mu(y')\geq k+1$; this is a contradiction.
\qed

\bigskip
\no
{\em Proof of Theorem \ref{covertheorem}.}
Let $V\sub\im(f)\sub \R^n$ be the set where the stable multiplicity 
function  $\mu$
is locally maximal; clearly $V$ is dense in $\im(f)$.
By Lemma \ref{semicontinuity}, $V$ is an open 
subset of $\R^n$, and $\mu$ is locally constant
on $V$. By Lemma \ref{allgreen}, the map $y\mapsto
\# f^{-1}(y)$ is a locally constant function on 
$V$.   It is therefore clear by Lemma
\ref{semicontinuity} that $f$ is locally
injective near any $x\in U:= f^{-1}(V)$, and hence
$f\restr_{U}$ is a covering map. If  
$W$ is a nonempty open set in $X$, then $f(W)$ has nonempty interior
by Proposition \ref{stablevalue}.  Hence $f(W)$ meets $V$, since $V$ is 
dense in $\im(f)$. It follows that $W$ meets $U=f^{-1}(V)$.  
This implies that $U$ is dense in $X$. 
\qed

\section{Weak Tangents}

In this section we briefly review some results on  weak tangents.
For more details see \cite{davsem} and \cite{BBI}. 

A {pointed metric space} is a pair $(Z, p)$, where $Z$ is a metric space
(with metric $d_Z$) and $p\in Z$.
A sequence  $(Z_k, p_k)$ of pointed metric spaces is said
to converge to a pointed metric space $(Z,p)$, if for every $R>0$ and for 
every $\eps>0$ there exist $N\in \N$,   a subset  $M\sub B_Z(p, R)$,
 subsets $M_k \sub B_{Z_k}(R)$ and bijections $f_k\: M_k\to M$ 
such that for $k\ge N$ 
\begin{itemize}
\item[(i)] $p\in M$,  $p_k\in M_k$, and  $f_k(p_k)=p$, 
\item[(ii)] the set $M$ is  $\eps$-dense in $B_{Z}(p, R)$, and 
the sets $M_k$ are  $\eps$-dense in $B_{Z_k}(p_k, R)$, 
\item[(iii)] $|d_{Z_k}(x,y)-d_Z(f_k(x), f_k(y))|<\eps$ whenever
$x,y\in M_k$. 
\end{itemize}

The definitions for pointed space convergence
 given in \cite{davsem} and \cite{BBI} are different, but equivalent.   

A complete 
metric space $S$  is called a {\em weak tangent} of the metric space
$Z$, if there exist a sequence  of  numbers 
$\la_k>0$  with $\la_k\to \infty$
for $k\to \infty$ and points $q\in S$, $p_k\in Z$
such that  the sequence of pointed spaces
 $(\la_kZ, p_k)$ converges to the pointed space $(S,q)$. Here we denote by 
$\la Z$ for $\la>0$ the metric space $(Z, \la d_Z)$. In other words,
$\la Z$ agrees with $Z$ as a set, but is equipped with  the  metric obtained
by rescaling the original metric  by the factor $\la>0$.
The set of all weak tangents of a metric space $Z$ is denoted by $\WT(Z)$. 
 If $X$, $Y$, $Z$ are metric spaces,   and $X$ is a weak tangent of 
$Y$ and $Y$ is a weak tangent of $Z$, then $X$ is a weak tangent
of $Z$, i.e., $X\in \WT(Y)$ and $Y\in \WT(Z)$ imply $X\in \WT(Z)$. 

A metric space $Z$ is called {\em uniformly perfect} if there exists
a constant $\la\ge 1$ such that for every $z\in Z$ and $0<R\le \diam(Z)$ 
we have $\bar B(z,R)\setminus B(z, R/\la)\ne \emptyset$.

For $Q>0$ we denote by ${\mathcal H}^Q$ the $Q$-dimensional Hausdorff measure
on a metric space $Z$. A complete metric  space $Z$ of positive diameter
 is called 
{\em Ahlfors $Q$-regular}, where  $Q>0$, if there exists 
a constant $C\ge 1$ such that 
$$  \frac{1}{C}R^Q  \le {\mathcal H}^Q (B(z,R)) \le CR^Q,$$
whenever $z\in Z$ and $0<R\le \diam(Z)$.     

A metric space $Z$ is called {\em doubling}, if there exists a number 
$N\in \N$ such that
every open  ball of radius $R$ in $Z$ can be covered by at 
most $N$ open balls of radius $R/2$. The space $Z$ is called {\em proper},
if closed balls in $Z$ are compact.  

Every Ahlfors regular space is uniformly perfect and doubling. 
A complete doubling space is proper. 
If $Z$ is a compact metric space  that is uniformly 
perfect and doubling, and $X\in \WT(Z)$, then $X$ is 
an unbounded  doubling  metric space. Since  $X$ is also  complete 
by definition, this space will be proper.

Suppose $f\: X\to Y$ is a map between a metric space $X$ and a 
doubling metric space  $Y$. The map
is called {\em regular} if it is Lipschitz and there exists 
a constant $N\in \N$ such that the inverse image of every open  ball $B$
 in $Y$  can be covered  by at most $N$ open balls   in $X$ 
with the same radius as $B$.

Note that this last condition implies that $f$ is of bounded multiplicity.
Indeed, we have $\#f^{-1}(y)\le N$ for $y\in Y$. For suppose that there
are $N+1$ distinct points $x_1, \dots, x_{N+1}\in f^{-1}(y)$.
Let $\eps>0$ be the minimum of the  distances $d_X(x_i, x_j)$ for $i\ne j$.
Consider the ball $B=B(y, \eps/2)$. By our assumption on  $f$ 
the preimage $f^{-1}(B)\supseteq f^{-1}(y)$ can be covered 
by $N$ open balls $B_1, \dots, B_N \sub X$ of radius $\eps/2$. But this is 
impossible, because each ball $B_i$ can contain at most one of the    points 
$x_1, \dots, x_{N+1}$.

The proof of the following proposition can be found in
\cite[Prop.\ 12.8]{davsem}.

\begin{proposition}\label{easydavsem}
 Let $X$ and $Y$ be  metric spaces, and $f\: X\to Y$
be a Lipschitz map. Suppose that $X$ is compact and
 Ahlfors $Q$-regular, where $Q>0$,
$Y$ is 
complete and doubling, and ${\mathcal H}^Q(f(X))>0$. 

Then there exist weak tangents $S\in \WT(X)$, $T\in \WT(Y)$, 
and a regular map $g\: S\to T$.
\end{proposition}

We will  need the following lemmas.

\begin{lemma}\label{biLip} Suppose   $X$ is a  metric space, 
and $f\:X\to \R^n$ is  regular. Assume that there is an 
open ball $B\sub \R^n$ and a set $U\sub f^{-1}(B)$ such that 
the map $g:=f\restr_U\: U\to B$ is a homeomorphism.
Then $g$ is a  bi-Lipschitz map. 
\end{lemma} 

It is understood that $U$ is equipped with the restriction of the metric
$d_X$ to $U$, and $B$ with  the Euclidean metric.

\proof Since $f$ is Lipschitz, the map $g$ is also Lipschitz. 
It remains to  obtain an upper bound for $d_X(x,y)$ in terms of  
 $|f(x)-f(y)|$,  whenever  
$x,y\in U$, $x\ne y$. Let $R:=2|f(x)-f(y)|>0$, $B':=B(x, R)$ and $S\sub 
B'\cap B$ be the Euclidean line segment connecting $f(x)$ and $f(y)$.
Then $E:=g^{-1}(S)$ is a compact connected set in $U$ containing 
$x$ and $y$. On the other hand, $E\sub f^{-1}(B')$.
If $N\in \N$ is associated with $f$ as in the definition of a regular map, 
then it follows that $E$ can be covered by $N$ open balls of radius 
$R$. Now we invoke the following elementary fact whose proof is left to the 
reader: If $E$ is a compact connected  set
in a metric space covered by open balls,
then  the diameter of $E$ is at most twice the sum of the radii of the balls.

In our situation we get the estimate
$$ d_X(x,y)\le \diam(E)\le 2N R= 4N |f(x)-f(y)|, $$
which proves that $g$ is a bi-Lipschitz homeomorphism.  
\qed

\begin{lemma} \label{WT1} Suppose $X$ and $Y$ are complete
 doubling metric spaces.
Suppose there exists a point $x\in X$, a neighborhood $U$ of
$x$ and a bi-Lipschitz map $f\: U\to V:=f(U)$ such that
$V$ is a neighborhood of $y:=f(x)$.
 
Then there exist $S\in \WT(X)$, $T\in \WT(Y)$, and a bi-Lipschitz
homeomorphism $g\: S\to T$.
\end{lemma}
 
The lemma says that under the given hypotheses the spaces
$X$ and $Y$ have bi-Lipschitz equivalent weak tangents.

\proof For $\la>0$ consider the pointed metric spaces $(\la U, x)$ 
and $(\la V, y)$, where $\la U$ and $\la V$ denote the metric spaces
whose underlying sets are $U$ and $V$ equipped with 
the restrictions  of the metric $d_X$ and $d_Y$, respectively, rescaled
by the factor $\la>0$. 
The map $f$ considered as a map between $(\la U, x)$ and $(\la V, y)$
preserves base points and  is bi-Lipschitz with a constant independent
of $\la$. Since $X$ and $Y$ are complete
and  doubling, it follows that in the terminology of David and Semmes
\cite[Sect.\ 8.5]{davsem} the mapping packages $f\: (\la U, x) \to (\la V, y)$
subconverge for $\la \to \infty$ to a mapping $g\: S\to T$. Here 
$S$  and $T$ are  limits  of the pointed spaces $(\la_k U, x)$
and $ (\la_k V, y)$, respectively, where $\la_k$ is a sequence of 
positive numbers with $\la_k\to \infty$ as $k\to \infty$. 
Since $U$ and $V$ are neighborhoods of $x$ and $y$, respectively,
 it follows that $S\in
\WT(X)$ and $Y\in \WT(Y)$ (cf.\ \cite[Lem.\ 9.12]{davsem}). Moreover, since the 
bi-Lipschitz constant of  $f\: (\la U, x) \to (\la V, y)$ is independent
of $\la$, the map $g$ will be bi-Lipschitz. 
There is a slight problem here, because it is not clear whether 
$g$ will be surjective. This problem can be 
addressed similarly as in the second part of the proof of 
Theorem \ref{TBA}. We may assume that the sequence $\la_k$ is 
such that not only the mapping packages $f\: (\la_k U, x) \to (\la_k V, y)$
converge, but also the mapping packages $f^{-1}\:
(\la_k V, y) \to (\la_k U, x)$,  to  $h\: T\to S$, say.
Then $g\circ h=\id_T$ which implies that $g$ is onto, and 
hence a bi-Lipschitz homeomorphism.  \qed

\section{Weak tangents and quasi-M\"obius actions}
\label{weaktang}

In this section we study weak tangents of compact metric spaces
which admit a uniformly \qm\ action for which the induced
action $G\acts\trip(Z)$ is cocompact.
As the reader will notice, 
all the results in this section remain true under the weaker 
assumption that every triple of distinct points in $Z$ can be blown
up to a uniformly separated triple  by a  uniform  \qm\ homeomorphism 
of $Z$, i.e., an $\eta$-\qm\ homeomorphism with $\eta$ independent of the
triple.

\begin{lemma}
\label{fillsin}
Suppose $Z$ is  a   uniformly perfect compact metric space,
 and $G\acts Z$  is an $\eta$-\qm\ action. 

\begin{itemize}
\item[(i)]
Suppose that for each  $k\in \N$ we are given  a set 
$D_k$ in a ball $B_k=B(p_k, R_k)\sub Z$ that is   $(\eps_k R_k)$-dense in 
$B_k$, where $\eps_k>0$, 
 distinct points $x_k^1, x_k^2, x_3^k \in 
B(p_k, \la_k R_k)$,  where $\la_k>0$, 
with
$$ d_Z(x_k^i, x_k^j)> \delta_k R_k
\for i,j\in \{1,2,3\},\,  i\ne j, $$
where $\delta_k>0$, 
and group elements $g_k\in G$ such that for $y^i_k:=g_k(x_k^i)$
we have 
$$d_Z(y_k^i, y_k^j)>
\delta' \for i,j\in \{1,2,3\},\, i\ne j, $$ 
where $\delta'>0$ is independent of $k$. 

Let $D'_k:=g_k(D_k)$, and suppose $\la_k\to 0$ for  $k\to \infty$
and that the sequence   $(\eps_k/\de^2_k)$ is bounded. 
Then 
$$ \dist_H(D'_k, Z)\to 0 \for k\to \infty. $$ 

\item[(ii)] 
Suppose in addition
 that $G\acts \trip(Z)$ is cocompact. 
 If $U\sub Z$ is a nonempty open set, then there 
exists a sequence $(g_k)$ in $G$ such that
$$\Diam(Z\setminus g_k(U)) \to 0 \for k\to \infty.  $$ 
\end{itemize}
\end{lemma}

In plain words (i)  essentially says that if we blow up
 a triple $(x^1,x^2,x^3)$ that lies in  a ball $B$
 to a uniformly separated triple, then a set $D$  in $B$ will be blown up 
to a  rather dense set in $Z$, if the triple $(x^1,x^2,x^3)$ lies deep 
inside $B$ and its    separation
is much larger than  $\dist_H(D,B)$. 

\medskip\noindent 
{\em Proof of} (i). Let $d=d_Z$. 
Consider fixed $k\in \N$ and
drop the subscript $k$ for simplicity. 
The image of a point $z\in Z$ under  $g=g_k$ will be denoted by
$z':=g(z)$.
Pick an arbitrary point in $Z$, and write
it in the form $x'=g(x)$ where $x\in Z$.  
We have to find a point in $D'$ close to $x'$.  

\smallskip\noindent
{\em Case 1: $x\in B(p,R)$.}
There is a point $y\in D\cap B$ with $d(x,y)\leq \eps R$.
Since the minimal distance between the points $x_1,x_2,x_3$
is at least $\de R$, we can find two of them, call them 
 $a$ and $b$, so that
$d(y,a)\ge \de R/2$ and
 $d(x,b)\ge \de R/2$. Hence 
$$\frac{d(x',y')d(a', b')}
{d( x',b')d( a', y')}\leq
\eta\left (\frac{d(x,y)d(a,b)}{d(x,b)d(a,y)}\right )
\le \eta(8\eps\la/\delta^2).$$
Rearranging factors, this implies that
$$ d(x', y')\leq \Diam(Z)^2 \eta(8\eps\la/\delta^2)/\delta'
\le C_1\eta(C_2\la).$$
The last expression becomes  uniformly small as $\la \to 0$. 

\smallskip\noindent
{\em Case 2: $x\not\in B(p,R)$.}
Since $\eps \lesssim \delta^2 \lesssim \la^2$,
we may assume that $\eps>0$ is small.  Then by the 
uniform perfectness of $Z$ and the $(\eps R)$-density 
of $D$ in $B$,  we can find a point $y\in D\cap B$
so that $d(y,p)/R$ is uniformly bounded away from zero,
$d(y,p)/R\ge c_0>0$  say. Note that $c_0$ does not depend on $k$.   
We may assume that $\la< c_0/2\le 1/2$. Then setting $a=x^1$ and 
$b=x^2$ we get   
\begin{eqnarray*}
\frac{d(x', y')d( a',b')}
{d(x', b' )d(a', y')}&\leq& 
\eta\left ( \frac{d(x,y)d(a,b)}{d(x,b)d(a,y)}\right) \\
&\leq& \eta\left(\frac{4\la d(x,p)}{(d(x,p)-\la R)(c_0-\la)}\right)\\
&\leq& \eta(16 \la/c_0).
\end{eqnarray*}
Rearranging factors, this implies that
$$d(x', y')\leq \Diam(Z)^2 \eta( 16 \la/c_0)/\delta'\le
C_3\eta(C_4\la).$$
Again the last expression becomes  uniformly small as $\la \to 0$. 

Since $y'\in D'$, the first part of the lemma follows. 

\medskip \noindent
{\em Proof of} (ii). Let $B=B(p, R)$ be a ball in 
$U$ with small radius $R\in (0,1/2]$. By the uniform 
perfectness of $Z$ we can find a triple $(x^1,x^2,x^3)$  of distinct points
in  $B(p,R^2)$ whose separation is comparable to $R^2$. 
Now use the cocompactness of $G\acts \Trip(Z)$ to find $g\in G$ mapping 
$(x^1,x^2,x^3)$ to a uniformly separated triple. 

Arguing as in Case 2 above, we find that whenever $x'$ and $
y'$ are points in $Z\setminus g(B(p, R))$, then 
$ d_Z(x', y')\lesssim \eta (CR).$ 
Hence $\diam (Z\setminus g(U))\lesssim \eta (C R)$, and the claim 
follows by making $R$ arbitrarily small.
\qed
\medskip

Before we state the next lemma we recall that in 
Section \ref{qm} we have defined a metric 
$\hat d_p$ on the one-point  compactification $\hat X$ of an unbounded 
locally compact pointed metric space $(X,p)$
associated with the metric $d=d_X$ and the base point $p$.

\begin{lemma}
\label{tangentonepoint}
Suppose $Z$ is  a
 compact metric space that is uniformly perfect and doubling, 
and $G\acts Z$  is a
uniformly \qm\ action for which the induced
action $G\acts\trip(Z)$ is cocompact.
 
If $(S,p)\in \WT(Z)$, then there exist a \qm\   homeomorphism 
$h\: (\hat S, \hat d_p) \to Z$.
Moreover, 
$h\restr_S\: S \to Z\setminus\{h(\infty)\}$ is also a    \qm\ homeomorphism.  
\end{lemma}
In other words, 
up to \qm\ homeomorphism  the space $Z$ 
is equivalent to the one-point compactification $\hat S$ 
of a weak tangent $(S,p)$ of $Z$ if we equip $\hat S$
with the canonical metric $\hat d_p$.
Conversely, up to \qm\ homeomorphism any  weak tangent of $Z$ is  
equivalent to $Z$ with one point removed.

\proof
Note that as a weak tangent of a uniformly perfect doubling metric space, 
 $S$ is unbounded and proper. 

From the definition of pointed space convergence it follows 
that for $k\in \N$  there exist   subsets
$\tilde D_k \sub B_S(p, k) \sub
S$ that are $(1/k)$-dense  in $B_S(p, k)$,
 numbers $\la_k>0$ with $\la_k\to \infty$, 
points $p_k\in Z$,  sets $D_k \sub B_{\la_k Z}(p_k, k) \sub \la_kZ$ that are 
$(1/k)$-dense in $B_{\la_k Z}(p_k, k)$
with respect to the metric $d_{\la_k Z}=\la_k d_Z$  
and bijections $f_k\: \tilde D_k\to  D_k$ such  that 
\begin{equation} \label{weaklimbilip} 
\frac12  d_S(x,y) \le \la_k d_Z(f_k(x), f_k(y))\le 2 d_S(x,y)
\for x,y\in \tilde D_k.
\end{equation}
Moreover, it  can be arranged  that 
each  set $\tilde D_k$ contains the points of a fixed triple $(q_1,q_2, q_3)\in 
\trip(S)$.  

Let $x_k^i:= f_k(q_i)$ for $i\in \{1,2,3\}$ and $k\in \N$.
Since the action  $G\acts \Trip(Z)$ is cocompact, for  $k\in \N$ we can find 
$g_k\in G$ such that the triples 
 $$(y^1_k, y^2_k,y^3_k):=g_k(x^1_k, x^2_k, x^3_k)
\in \trip(Z)$$
are uniformly separated.

The density condition for the sets $D_k$
rephrased in terms of the metric 
$d_Z$ says that $ D_k$ is 
$(\la_k/k)$-dense in $B_Z(p_k, \la_k k)$  with respect to $d_Z$.    
Moreover, in terms of the metric $d_Z$, 
 the triple $(x^1_k, x^2_k, x^3_k)$ has separation comparable
to $\la_k$ and is contained in a ball  centered at $p_k$ whose radius is also
comparable to $\la_k$.   
It follows from Lemma \ref{fillsin}  that  for $D'_k:=g_k(D_k)$ we have 
\begin{equation}\label{ddk1}
 \lim_{k\to \infty} \dist_H(D'_k, Z)=0, 
\end{equation}
where $\dist_H$ refers to the Hausdorff distance in $Z$.

The density condition  for the sets $\tilde D_k \sub S \sub \hat S $
and the inequality \eqref{metricAC} 
for the metric $\hat d_p$ imply that 
\begin{equation}\label{ddk2}
 \lim_{k\to \infty} \dist_H(\tilde D_k, \hat S)=0,
\end{equation}
where $\dist_H$ refers to the Hausdorff distance in $(\hat S, \hat d_p)$.

Consider the maps $h_k\: (\tilde D_k, \hat d_p\restr_{\tilde D_k})\to Z$
defined by 
$h_k(x)=g_k(f_k(x))$ for $x\in \tilde D_k$. Note that 
it follows from 
Lemma \ref{dhat}, inequality \eqref{weaklimbilip}
and the fact that the action $G\acts 
Z$ is uniformly \qm\ that the maps $h_k$  are $\eta$-\qm\ with $\eta$
independent of $k$. 
Moreover, each map $h_k$ maps the triple $(q_1,q_2,q_3)$ to 
the uniformly separated triple $(y^1_k, y^2_k,y^3_k)$. 
Finally, $D'_k=h_k(\tilde D_k)$ and so by \eqref{ddk1} and \eqref{ddk2}
we can  apply
Lemma~\ref{TBA}. It follows  that the sequence $(h_k)$ subconverges 
to a \qm\ homeomorphism $h\: (\hat S, \hat d_p)\to Z$. 

The second part of the lemma follows by
observing  $h\restr_S\: S \to Z\setminus\{f(\infty)\}$ is \qm, 
since this map 
 the composition 
of the maps $\id_S\: S \to (S, \hat d_p\restr_S)$ which is \qm\ by 
Lemma~\ref{dhat}  and the map  
$h\restr_S\: (S, \hat d_p\restr_S) \to Z\setminus \{h(\infty)\}$ which is 
\qm\ by the first part of the proof.  
\qed

\begin{lemma}
\label{dim}
Suppose $Z$ is  a
 compact  metric space that is uniformly perfect and doubling,
and $G\acts Z$  is a
uniformly \qm\ action for which the induced
action $G\acts\trip(Z)$ is cocompact.

If $\topdim(Z)=n\in \N$, then
$\topdim(U)=n$ whenever $U$ is a  nonempty open subset  of $Z$ 
or of any weak tangent of $Z$.
\end{lemma}
\proof
If $U\sub Z$ is a nonempty open set, we can find a 
nonempty open set $V$ with $\bar V\sub U$.
By Lemma \ref{fillsin} there is a sequence $(g_k)$ in $G$  such
that $\diam(Z\setminus g_k(\bar V)) \to 0$ for 
$k\to \infty$.  Hence the complement of $\bigcup_{k\in N} g_k(\bar V)$ in $Z$ 
can contain at most one point. 
 Topological dimension is  invariant under homeomorphisms, and  
and  does not 
increase under a  countable union of closed sets
(cf.\ \cite[Thm.\ II.\ 1]{nag}). So  we get
$\topdim(Z)\leq \topdim (\bar V)\leq \topdim(U)\leq  \topdim(Z)$.

If $U$ is a nonempty open subset 
of any weak tangent $S$  of $Z$, then $U$ is also 
 an   open subset of the  one-point compactification
of $S$.
Hence by Lemma \ref{tangentonepoint}, the set  $U$ is homeomorphic to 
a nonempty open subset of $Z$. Therefore 
$\topdim(U)=\topdim(Z)$
by the first part of the proof. 
\qed

\begin{lemma}
\label{rescaling}
Suppose $X$ and $Y$ are  compact metric spaces that are uniformly 
perfect and doubling, and suppose   
 $G\acts Z$ and $H\acts X$ are  uniformly \qm\
actions for which the induced actions $G\acts\trip(X)$,
$H\acts\trip(Y)$ are cocompact.

If there exist $S\in \WT(X)$ and $T\in \WT(Y)$ and a \qs\
homeomorphism $f\: S\to T$, then there exists a \qm\ homeomorphism
$g\: X\to Y$.
\end{lemma}

So if $X$ and $Y$ have weak tangents that are quasi-symmetrically equivalent, 
then $X$ and $Y$ are equivalent up to a  \qm\ homeomorphism.

\proof  Let $p$ and $q$ be the base points in $S$ and $T$, respectively, 
and consider the one-point compactifications $(\hat S, \hat d_p)$ 
and $(\hat T, \hat d_q)$.  If we define $\hat f(x)=f(x)$ for 
$x\in X$, and $\hat f(\infty)=\infty$, then 
 \eqref{metricAC} implies that  $\hat f\: (\hat S, \hat d_p) \to
(\hat T, \hat d_q)$ is a  \qm\ homeomorphism.
 Since $(\hat S, \hat d_p)$ is equivalent
to $X$  and 
$(\hat T, \hat d_q)$ is  equivalent to $Y$ up to 
\qm\ homeomorphisms by Lemma~\ref{tangentonepoint}, the claim follows.
\qed

\section{Proof of Theorem \ref{mainthm}}
\label{proofmainthm}

Let $Z$ and  $G\acts Z$ be as in the statement of Theorem
\ref{mainthm}.  

We are given that $\topdim(Z)=n$.  This implies \cite[Thm.\ III. 1]{nag}
that there is a continuous  map $f_0\:Z\ra\S^n$ with
a stable value $y\in \S^n$; in fact any continuous map $f_1\: Z\ra\S^n$
for which $\dist(f_0,f_1)$ is sufficiently small 
will also have $y$ as a stable value.

Every continuous function $g_0\: Z \to \R$ can be approximated by a 
Lipschitz function $g_1\: Z\to \R$ such that $\dist(g_0, g_1)$ is arbitrarily
small. This standard fact can be established by using  Lipschitz 
partitions  of unity in $Z$ subordinate to a cover of $Z$ by small 
balls with controlled overlap.
We apply this  to the $n+1$ coordinate functions
of the map $f_0\:Z\ra\S^n \sub \R^{n+1}$  to obtain  Lipschitz maps  
   on $Z$ which are arbitrarily 
close to $f_0$ and map $Z$ into small neighborhoods 
of $\S^n$ in  $\R^{n+1}$. Composing these  maps  with the radial projection 
from the origin in  $\R^{n+1}$ to $\S^n$, we can find 
 Lipschitz maps from $Z$ into $\S^n$  arbitrarily  
close to $f_0$. In particular, there exists 
a Lipschitz map $f\: Z\to \S^n$ such  that $y$ is a stable value of 
$f$. Then $\im(f)$ is a neighborhood of $y$, and so 
${\mathcal H}^n(\im(f))>0$.

We now apply Proposition~\ref{easydavsem}  to obtain 
a weak tangent $S$ of $Z$, a weak tangent $T$ of $\S^n$  and 
a regular
map $\phi\:S \to T$. Note that every weak tangent of $\S^n$ is 
isometric to $\R^n$, and so $T=\R^n$.

As we have seen, the fact that $\phi$ 
is regular implies  that 
$\phi$ has bounded multiplicity.
By Lemma \ref{dim}, every nonempty open subset of $S$ has
topological dimension $n$.   Therefore, 
by Theorem \ref{covertheorem} (applied 
 to the closure of the some 
bounded nonempty open  set in $S$ as the space $X$)
 there is a nonempty open subset 
$U\sub S$ such that $\psi:=\phi\restr_U$ is a homeomorphism
onto an open subset of $\R^n$. Shrinking the open set $U$ if necessary, 
we may assume that $\psi$ is a homeomorphism onto an open ball 
$B$ in $\R^n$.
Now Lemma~\ref{biLip} shows that $\psi$ is bi-Lipschitz. 
Choosing $x\in U$ and setting $y:=\psi(x)$ we are in the situation of
Lemma~\ref{WT1}.
We conclude that $S$  has a weak tangent 
bi-Lipschitz equivalent to a  weak tangent of $\R^n$. 
 Since the weak tangents
of $S$ are also weak tangents of $Z$, and all weak tangents of $\R^n$ 
are isometric to $\R^n$,  we  see that
$Z$ and $\S^n$ have bi-Lipschitz equivalent weak tangents.
Since the group of M\"obius transformations  induces a uniformly 
\qm\ action on $\S^n$ and a cocompact action on $\trip(\S^n)$, 
Lemma~\ref{rescaling} implies that 
there exists
 a \qm\ homeomorphism $h\:Z\ra\S^n$.
As a \qm\ homeomorphism between bounded spaces, the map $h$ will also be 
quasi-symmetric. 
 Conjugating
the uniformly \qm\ action $G\acts Z$ by $h$, we get
a  uniformly \qm\  
action $G\acts \S^n$ such that the induced action 
$G\acts \trip(\S^n)$ is cocompact. 
By  a result of Tukia~\cite[Cor.\ G(a)]{Tuk}, this action
is conjugate by a quasiconformal homeomorphism to
an action by M\"obius transformations.   Since quasiconformal
homeomorphisms of $\S^n$ onto itself are \qs,  Theorem \ref{mainthm}
follows.
\qed

\medskip

The method of proving Theorem~\ref{mainthm} also  leads to the following 
result.    

\begin{theorem}
\label{answer}
Let $n\in\N$, and let $Z$ be a compact,
Ahlfors $n$-regular metric space
of topological dimension $n$. Suppose every triple of distinct 
points in $Z$ can be mapped to a uniformly separated triple by 
a uniform \qm\  homeomorphism of $Z$.
Then $Z$ is quasi-symmetrically equivalent to the standard sphere $\S^n$.
\end{theorem}

\proof In the same way as in the proof of Theorem~\ref{mainthm}, we 
see that $Z$ has a weak tangent  bi-Lipschitz equivalent to $\R^n$. 
As we remarked in the beginning of Section~\ref{weaktang}, the results
in this section remain true if the assumption on the group 
action is replaced by  the assumption 
that every triple of distinct
points in the space under consideration 
  can be mapped to a uniformly separated triple by
a uniform \qm\  homeomorphism. So by the analog of Lemma~\ref{rescaling},
 we again  
obtain  a \qm,  and hence  
 quasi-symmetric, 
homeomorphism $h\: Z \ra \S^n$. \qed 

\medskip 

This theorem justifies the remark in the introduction about
the question of Heinonen and Semmes---recall that \qm\ homeomorphisms 
of compact metric spaces  are quasi-symmetric.   We see that the
three point homogeneity condition can be relaxed to a ``cocompact
on triples'' condition, at the cost of requiring the homeomorphisms
to be uniformly  \qm.

\section{$\CAT(-1)$-spaces and isometric  group actions}
\label{hypgp}

We refer the reader to \cite{ghysdela} for general
background on Gromov  hyperbolic spaces.
 
A metric space $X$ is called {\em geodesic},
if any two points $x,y\in X$ can be 
joined by a geodesic segment in $X$, i.e., a curve whose length is equal to 
the distance of $x$ and $y$. In the following we will always assume 
that $X$ is proper and geodesic. 

Let  $X$ be a Gromov hyperbolic space, and $\geo X$ be 
its boundary at infinity. There is a natural topology on $X\cup \geo X$
making this union compact. 
If $p\in X$, $a,b\in\geo X$, we let
$[a,b]_p$ denote the Gromov product of $a,b\in\geo X$
with respect to the base point $p$.  When $c>0$ is sufficiently
small,  the function
\begin{equation}
\label{boundarydistance}
d(a,b)\defeq \exp (-c[a,b]_p)
\end{equation}
is equivalent up to a multiplicative factor to a metric
on $\geo X$; any two metrics of this type are \qs ally
equivalent by the identity map.   
Fix one such metric on $\geo X$. If we denote the group of 
isometries of $X$ by $\Isom(X)$, then 
we get an induced  action $\Isom(X)\acts\geo X$
which is a uniformly \qm\ action,
 \cite[Prop.\ 4.5]{Pau}.  In fact, every
quasi-isometry $f\:X\ra X$  induces an $\eta$-\qm\
homeomorphism $\geo X\ra\geo X$ where $\eta$ depends only
the parameters of the quasi-isometry and the hyperbolicity constant 
of $X$.
 
Now suppose that $X$ is a $\CAT(-1)$-space (see \cite{Bou1} for 
more details on the topics discussed in the following).
Then  for every $p\in X$
we get a canonical metric on $\geo X$ as follows.
For every point  $a\in \geo X$, there is a unique 
geodesic ray $\ol{pa}$ starting at $p$ whose asymptotic class
represents $a$.
  Let  
 $a,b\in \geo X$, and consider points
$x\in\ol{pa}$, $y\in\ol{pb}$.
  Let $\De \tilde{p}\tilde{x}\tilde{y}$
be a comparison
triangle (in the hyperbolic plane) for the triangle
$\De pxy$, and let $\cangle_p(x,y)$ denote the
angle at $\tilde{p}$.   When $x$ and $y$ tend to
infinity along the rays $\ol{pa}$ and $\ol{pb}$, respectively, 
 the comparison angle $\cangle_p(x,y)$
has a limit, which we define to be the distance
between $a$ and $b$. 
This metric agrees up to a bounded factor with
the expression in (\ref{boundarydistance}) when
$c=1$.

Suppose $G\acts X$ is an isometric action of a group on 
a $\CAT(-1)$-space $X$. If $x\in X$, then we denote its orbit 
under $G$ by 
$$Gx:=\{g(x): g\in G\}.$$ 
The {\em limit set}  $\La(G) \sub  \geo X$ of $G$ is by definition 
the set of all accumulation points of an orbit $Gx$ on $\geo X$. 
This set is independent of $x\in X$. \
The group action $G\acts X$ is called {\em properly discontinuous} 
if 
$$ \{ g\in G: g(K) \cap K\ne \emptyset\}$$   
is finite for every compact subset $K$ of $X$.

A subset $Y\sub X$ is {\em quasi-convex} if there is  a
constant $C$ such that any geodesic segment with endpoints
in $Y$ lies in the $C$-neighborhood of $Y$.
The action $G\acts X$ is {\em quasi-convex cocompact} if there
is a $G$-invariant quasi-convex subset $Y\subseteq X$ on which 
$G$ acts with compact quotient $Y/G$. 
The group  $G$ is quasi-convex cocompact if and only if all  orbits 
$Gx$ are  quasi-convex.

We will need the following result due to Bourdon
\cite[0.3 Th\'eor\`eme ($\H^n$ case)]{Bou2}. 

\begin{theorem}
\label{bourdon}
 Let  $n\ge 2$, $G$ be  a group, and $X$  a 
$\CAT(-1)$-space. Suppose we have isometric group actions 
$G\acts X$ and $G\acts \H^{n+1}$ which are properly discontinuous.
Suppose that $G\acts X$ is quasi-convex 
cocompact and $G\acts\H^{n+1}$ is cocompact. 
If the Hausdorff dimension of $\La(G)\sub \geo X$ is equal to $n$,
then there exists a $G$-equivariant isometry of $\H^{n+1}$ onto 
a  convex, $G$-invariant set $Y\sub X$.
\end{theorem}

Actually, Bourdon proved this under the additional
 assumption that the group action 
$G\acts \H^{n+1}$ is faithful. In this case $G$ is isomorphic
to  a uniform lattice in $\Isom(\H^{n+1})$. The proof of the above more 
general
version is the same as the proof of his original result.  
 
\medskip\noindent
{\em Proof of Theorem \ref{cat-1theorem}.}
Consider the  induced actions   $G\acts \La(G)$ and 
$G\acts \Trip(\La(G))$. Since $G\acts X$ is isometric, 
$G\acts \La(G)$ is uniformly \qm.  
 Since 
the action $G\acts X$ is properly discontinuous, the same is true 
for $G\acts \Trip(\La(G))$.
Moreover, since $G\acts X$ is quasi-convex cocompact, $G\acts \Trip(\La(G))$
is  cocompact.  

 Since the Hausdorff dimension
of $\La(G)$ is $n$, this space will actually be
Ahlfors $n$-regular (cf.\ \cite[Section~7]{coo}).
Now $n$ is also the topological dimension of $\La(G)$ by assumption.
 By Theorem \ref{mainthm}, the action $G \acts \La(G)$ is
quasi-symmetrically  conjugate  to an 
action $G\acts \S^n$ by M\"obius transformations.
The action $G\acts \Trip(\S^n)$ is  properly
discontinuous and cocompact.
This implies that there is a  properly discontinuous, cocompact, and
isometric action $G\acts \H^{n+1}$ which induces
the action $G\acts \S^n=\geo \H^{n+1}$.
Since 
 $n\ge 2$ we can apply Bourdon's theorem, and 
conclude that there exists a $G$-equivariant 
isometric embedding of $\H^{n+1}$ onto a convex, $G$-invariant  set 
 $Y\sub X$ on which $G$ acts cocompactly. The result follows.   
\qed

\medskip
As the proof shows, $n\ge 2$ is only used in the last step. 
In particular, even in the case $n=1$ we can  still conclude that 
$\La(G)$ is quasi-symmetrically equivalent to $\S^1$, and 
that there is an action $G\acts \H^2$ which 
isometric, properly discontinuous and cocompact.

\bibliography{refs}

\begin{thebibliography}{10}

\bibitem{Bou1}
{\sc M.~Bourdon}, {\em Structure conforme au bord et flot g\'eod\'esique d'un
  {C}{A}{T}(-1)-espace}, Enseign. Math. (2), 41 (1995), pp.~63--102.

\bibitem{Bou2}
\leavevmode\vrule height 2pt depth -1.6pt width 23pt, {\em Sur le birapport au
  bord des {C}{A}{T}(-1)-espaces}, Inst. Hautes \'Etudes Sci. Publ. Math., 83
  (1996), pp.~95--104.

\bibitem{bowen}
{\sc R.~Bowen}, {\em Hausdorff dimension of quasicircles}, Inst. Hautes
  \'Etudes Sci. Publ. Math., 50 (1979), pp.~11--25.

\bibitem{BBI}
{\sc D.~Burago, Y.~Burago, and S.~Ivanov}, {\em A course in metric geometry},
  American Mathematical Society, Providence, RI, 2001.

\bibitem{coo}
{\sc M.~Coornaert}, {\em Mesures de {P}atterson-{S}ullivan sur le bord d'un
  espace hyperbolique au sens de {G}romov}, Pacific J. Math., 159 (1993),
  pp.~241--270.

\bibitem{davsem}
{\sc G.~David and S.~Semmes}, {\em Fractured fractals and broken dreams},
  Oxford University Press, New York, 1997.

\bibitem{ghysdela}
{\sc {\'E}.~Ghys and P.~de~la Harpe}, eds., {\em Sur les groupes hyperboliques
  d'apr\`es {M}ikhael {G}romov}, Birkh\"auser, Boston, MA, 1990.

\bibitem{heinsem}
{\sc J.~Heinonen and S.~Semmes}, {\em Thirty-three yes or no questions about
  mappings, measures, and metrics}, Conform. Geom. Dyn., 1 (1997), pp.~1--12
  (electronic).

\bibitem{nag}
{\sc J.~Nagata}, {\em Modern {D}imension {T}heory}, Heldermann Verlag, Berlin,
  1983.

\bibitem{Pau}
{\sc F.~Paulin}, {\em Un groupe hyperbolique est d\'etermin\'e par son bord},
  J. London Math. Soc. (2), 54 (1996), pp.~50--74.

\bibitem{sul1}
{\sc D.~Sullivan}, {\em Discrete conformal groups and measurable dynamics},
  Bull. Amer. Math. Soc. (N.S.), 6 (1982), pp.~57--73.

\bibitem{Tuk}
{\sc P.~Tukia}, {\em On quasiconformal groups}, J. Analyse Math., 46 (1986),
  pp.~318--346.

\bibitem{Vai}
{\sc J.~V{\"a}is{\"a}l{\"a}}, {\em Quasi-{M}\"obius maps}, J. Analyse Math., 44
  (1984/85), pp.~218--234.

\bibitem{yue}
{\sc C.~Yue}, {\em Dimension and rigidity of quasi-{F}uchsian representations},
  Ann. of Math. (2), 143 (1996), pp.~331--355.

\end{thebibliography}
\bibliographystyle{siam}
\addcontentsline{toc}{subsection}{References}

\end{document}